\documentclass[12pt, twoside, eqno]{article}
\usepackage{latexsym}
\usepackage{amssymb}
\usepackage{amsfonts}
\usepackage{amsmath}
\begin{document}

\noindent {\bf \large How we pass from semigroups to  
hypersemigroups}\bigskip

\noindent{\bf Niovi Kehayopulu}\bigskip

\noindent{\bf Abstract.} {\small The aim of writing this paper is 
given in the title. The results on semigroups can be easily 
transferred to hypersemigroups in the way indicated in the present 
paper.}\section{Introduction and prerequisites}An {\it hypergroupoid} 
is a nonempty set $H$ with an hyperoperation $$\circ : H\times H 
\rightarrow {\cal P}^*(H) \mid (a,b) \rightarrow a\circ b$$on $H$ and 
an operation $$* : {\cal P}^*(H)\times {\cal P}^*(H) \rightarrow 
{\cal P}^*(H) \mid (A,B) \rightarrow A*B$$ on ${\cal P}^*(H)$ 
(induced by the operation of $H$) such that 
$$A*B=\bigcup\limits_{(a,b) \in\,A\times B} {(a\circ b)}$$ for every 
$A,B\in {\cal P}^*(H)$ (${\cal P}^*(H)$ denotes the set of nonempty 
subsets of $H$). As the operation ``$*$" depends on the 
hyperoperation ``$\circ$", an hypergroupoid can be also denoted by 
$(H,\circ)$ (instead of $(H,\circ,*)$). If $H$ is an hypergroupoid 
then, for every $x,y\in H$, we have $\{x\}*\{y\}=x\circ y.$
A nonempty subset $A$ of an hypergroupoid $(H,\circ)$ is called a 
{\it left} (resp. {\it right}) {\it ideal} of $H$ if $H*A\subseteq A$ 
(resp. $A*H\subseteq A$). A subset of $H$ which is both a left ideal 
and a right ideal of $H$ is called an {\it ideal} of $H$. A nonempty 
subset $A$ of $H$ is called a {\it subsemigroup} of $H$ if 
$A*A\subseteq A$. Clearly, every left (resp. right) ideal of $H$ is a 
subsemigroup of $H$. An hypergroupoid $H$ is called {\it 
hypersemigroup} if $\{x\}*(y\circ z)=(x\circ y)*\{z\}$ for every 
$x,y,z\in H$. This is equivalent to saying that 
$\{x\}*{\Big(}\{y\}*\{z\}{\Big)}={\Big(}\{x\}*\{y\}{\Big)}*\{z\}$  
for every $x,y,z\in H$. Identifying the $\{x\}$ by $x$ and the 
$\{z\}$ by $z$ we could write, for short, $x*(y\circ z)=(x\circ y)*z$ 
for all $x,y,z\in H$. If $H$ is an hypersemigroup and $A,B,C\in {\cal 
P}^*(H)$, then we 
have\begin{eqnarray*}(A*B)*C&=&\bigcup\limits_{(a,b,c)\in A\times 
B\times C} {{\Big(}(a\circ b)*\{ c\}\Big) 
}\\&=&\bigcup\limits_{(a,b,c) \in A\times B\times 
C}{\Big(}\{a\}*(b\circ c){\Big)}=A*(B*C)\\&=&\bigcup\limits_{(a,b,c) 
\in A\times B\times 
C}{\Big(}\{a\}*\{b\}*\{c\}{\Big)}.\end{eqnarray*}Thus we can write  
$(A*B)*C=A*(B*C)=A*B*C$.\\From this, using induction, follows that 
for any finite family $A_1,A_2, ..., A_n$ of elements of ${\cal 
P}^*(H)$, we have$$A_1*A_2* ..... * 
A_n=\bigcup\limits_{({a_1},{a_2}...{a_n} )\in {A_1} \times {A_2} 
\times ... \times {A_n}}{\Big(} {\{ {a_1}\}*\{ {a_2}\}* ,... * \{ 
{a_n}} \}{\Big)}.$$The following proposition, though clear, plays an 
essential role in the theory of hypergroupoids.\medskip

\noindent{\bf Proposition 1.1.} {\it Let $(H,\circ)$ be an 
hypergroupoid, $x\in H$ and $A,B\in {\cal P}^*(H)$. Then we have the 
following:

$(1)$ $x\in A*B$ $\Longleftrightarrow$ $x\in a\circ b$ for some $a\in 
A$, $b\in B$.

$(2)$ If $a\in A$ and $b\in B$, then $a\circ b\subseteq A*B$.} 
\medskip

Following Zadeh, any mapping $f : H\rightarrow [0,1]$ of an 
hypergroupoid $H$ into the closed interval $[0,1]$ of real numbers is 
called a {\it fuzzy subset} of $H$ (or a {\it fuzzy set} in $H$) and 
$f_A$ (: the characteristic function of $A$) is the mapping$$f_A : H 
\rightarrow [0,1] \mid x \rightarrow f_A (x)=\left\{ 
\begin{array}{l}
1\,\,\,\,\,$if$\,\,\,\,x \in A\\
0\,\,\,\,$if$\,\,\,\,x \notin A.
\end{array} \right.$$The concepts of fuzzy left and fuzzy right 
ideals of semigroups can be naturally transferred in case of an 
hypergroupoid as follows: Let $H$ be an hypergroupoid. A fuzzy subset 
$f$ of $H$ is called a {\it fuzzy left ideal} of $H$ if
$f(x\circ y)\ge f(y)$ for all $x,y\in H$, in the sense that $x,y\in 
H$ and $u\in x\circ y$ implies $f(u)\ge f(y)$. A fuzzy subset $f$ of 
$H$ is called a {\it fuzzy right ideal} of $H$ if
$f(x\circ y)\ge f(x)$ for all $x,y\in H$, meaning that if $x,y\in H$ 
and $u\in x\circ y$, then $f(u)\ge f(x)$. A fuzzy subset $f$ of $H$ 
is called a {\it fuzzy ideal} of $H$ it is both a fuzzy left and a 
fuzzy right ideal of $H$. As one can easily see, a fuzzy subset $f$ 
of $H$ is a fuzzy ideal of $H$ if and only\ $f(x\circ y)\ge 
\max\{f(x),f(y)\} \mbox { for all } x,y\in H$, in the sense that 
$x,y\in H \mbox { and } u\in x\circ y \mbox { implies } f(u)\ge 
\max\{f(x),f(y)\}.$

Fuzzy interior ideals of semigroups (without order) and fuzzy simple 
semigroups (without order) have been first considered by Kuroki in 
[2]. If $S$ is a regular or an intra-regular ordered semigroup, then 
the ideals and the interior ideals of $S$ coincide. In particular, if 
$S$ is a $poe$-semigroup (i.e. an ordered semigroup having a greatest 
element with respect to the order), then the ideal elements and the 
interior ideal elements of $S$ coincide as well. In an attempt to 
show the similarity between the theory of ordered semigroups and the 
theory of fuzzy ordered semigroups, we have proved in [1] that in 
regular and in intra-regular ordered semigroups the fuzzy ideals and 
the fuzzy interior ideals coincide. We also proved that a set $A$ is 
an interior ideal of an ordered semigroup $S$ if and only its 
characteristic function $f_A$ is a fuzzy interior ideal of $S$. We 
finally proved that an ordered semigroup is simple if and only if it 
is fuzzy simple, equivalently, if every fuzzy interior ideal of $S$ 
is a constant function. The present paper is based on the papers [1] 
and [2] in the References, and its aim is to show the way we pass 
from semigroups to hypersemigroups. Analogous results hold for 
ordered hypersemigroups as well. Results on semigroups can be 
translated into the language of the hypersemigroups in the same way. 
Although we cannot  say ``all the results", up today, I have never 
seen a result that cannot be translated. And this is why, when we 
work on an hypersemigroup, we should indicate inside the paper the 
1--3 papers on semigroups on which our paper is based.  Instead, the 
corresponding papers on semigroups are sometimes (not always) cited 
in the References among many other papers on hypersemigroups (or even 
on semigroups) not related with the investigation. Concerning the 
present paper, as one can see, only the papers on semigroups cited in 
the References are needed in the investigation and nothing else. We 
use the terms left (right) ideal, bi-ideal, quasi-ideal instead of 
left (right) hyperideal, bi-hyperideal, quasi-hyperideal and so on, 
and this is because in this structure there are no two kind of left 
ideals, for example, to distinguish them as left ideal and left 
hyperideal. The left ideal in this structure is that one which 
corresponds to the left ideal of semigroups. We will give further 
interesting information about the hypersemigroups (without order) in 
a forthcoming paper.\newpage\section{Main results}{\bf Definition 
2.1.} Let $H$ be an hypersemigroup. A nonempty subset $A$ of $H$ is 
called an {\it interior ideal} of $H$ if$$H*A*H\subseteq A.$$By a 
{\it subidempotent interior ideal} of $H$ we mean an interior ideal 
of $H$ which is at the same time a subsemigroup of $H$. \medskip

\noindent{\bf Definition 2.2.} Let $H$ be an hypergroupoid. A 
nonempty subset $A$ of $H$ is called a {\it fuzzy interior ideal} of 
$H$ if

$$f{\Big(}(x\circ a)*\{y\}{\Big)}\ge f(a) \mbox { for every } 
x,a,y\in H.$$For an hypersemigroup, we clearly have
$$(x\circ a)*\{y\}=\{x\}*(a\circ y)=\{x\}*\{a\}*\{x\}.$$
{\bf Proposition 2.3.} {\it Let H be an hypersemigroup. If A is an 
interior ideal of H, then the characteristic function $f_A$ is a 
fuzzy interior ideal of H. ``Conversely", if $A$ is a nonempty subset 
of $H$ such that $f_A$ is a fuzzy interior ideal of $H$, then $A$ is 
an interior ideal of $H$}.\medskip

\noindent{\bf Proof.} $\Longrightarrow$. Let $x,a,y\in H$. Then
$f_A{\Big(}(x\circ a)*\{y\}{\Big)}\ge f_A(a)$. In fact: Let $u\in 
(x\circ a)*\{y\}$. If $a\in A$, then $f_A(a)=1$. Since $A$ is an 
interior ideal of $H$, we have $H*A*H\subseteq A$. So we have $u\in 
\{x\}* \{a\}* \{y\}\subseteq H*A*H\subseteq A$. Then $u\in A$, and 
$f_A(u)=1$. Thus we get $f_A(u)\ge f_A(a)$. Let now $a\notin A$. Then 
$f_A(a)=0$. Since $f_A$ is a fuzzy subset of $H$ and $u\in H$, we 
have $f_A(u)\ge 0$. Thus we have $f_A(u)\ge f_A(a)$.\\
$\Longleftarrow$. Let $A$ be a nonempty subset of $H$ and $f_A$ a 
fuzzy interior ideal of $H$. Then $H*A*H\subseteq A$. Indeed: Let 
$u\in H*A*H$. Then $u\in (x\circ a)*\{y\}$ for some $x,y\in H$, $a\in 
A$. Since $f_A$ a fuzzy interior ideal of $H$, we have $f_A(u)\ge 
f_A(a)=1$. Since $f_A$ is a fuzzy subset of $H$ and $u\in H$, we have 
$f_A(u)\le 1$. So we have $f_A(u)=1$, and $u\in A$. 
$\hfill\Box$\medskip

\noindent{\bf Proposition 2.4.} {\it Let $H$ be an hypersemigroup. If 
f is a fuzzy ideal of H, then f is a fuzzy interior ideal of 
H.}\medskip

\noindent{\bf Proof.} Let $x,a,y\in H$. Then $f{\Big(}(x\circ 
a)*\{y\}{\Big)}\ge f(a)$. In fact: Let $u\in (x\circ a)*\{y\}$. By 
Proposition 1.1, there exists $v\in x\circ a$ such that $u\in v\circ 
y$. Since $f$ is a fuzzy left ideal of $H$, we have $f(x\circ a)\ge 
f(a)$ and, since $v\in x\circ a$, we get $f(v)\ge f(a)$. Since $f$ is 
a fuzzy right ideal of $H$, we have $f(v\circ y)\ge f(v)$ and, since 
$u\in v\circ y$ we get $f(u)\ge f(v)$. Then we have $f(u)\ge f(a)$, 
and the proof is complete. $\hfill\Box$\medskip

\noindent{\bf Definition 2.5.} An hypersemigroup $H$ is called {\it 
regular} if for every $a\in H$ there exists $x\in H$ such that 
$a\in\{a\}*(x\circ a)$. \medskip

\noindent{\bf Proposition 2.6.} {\it Let H be an hypersemigroup. The 
following are equivalent:

$(1)$ H is regular.

$(1)$ $a\in \{a\}*\{x\}*\{a\}$ for every $a\in H$.

$(2)$ $A\subseteq H*A*H$ for every nonempty subset of $H$. }\medskip

\noindent{\bf Proposition 2.7.} {\it Let H be a regular 
hypersemigroup and A an interior ideal of H. Then A is a subsemigroup 
of H.}\medskip

\noindent{\bf Proof.} Since $A$ is an interior ideal of $H$, we have 
$H*A*H\subseteq A$. Since $H$ is regular, we have $A\subseteq A*H*A$. 
Then we have 
\begin{eqnarray*}A*A&\subseteq&(A*H*A)*A=(A*H)*A*A\subseteq 
(H*H)*A*H\\&\subseteq&H*A*H\subseteq A,\end{eqnarray*}so $A$ is a 
subsemigroup of $H$. $\hfill\Box$\medskip

\noindent{\bf Proposition 2.8.} {\it Let H be a regular 
hypersemigroup and f a fuzzy interior ideal of H. Then f is a fuzzy 
ideal of H}.\medskip

\noindent{\bf Proof.} Let $a,b\in H$ $(\Rightarrow f(a\circ b)\ge 
f(a) \mbox { and } f(a\circ b)\ge f(b) \;\;?\;)$\\
Let $u\in a\circ b$. Then $f(u)\ge f(a)$. Indeed: Since $a\in H$ and 
$H$ is regular, there exists $x\in H$ such that 
$a\in\{a\}*\{x\}*\{a\}$. Then$$a\circ b\subseteq 
\{a\}*\{x\}*\{a\}*\{b\}=(a\circ x)*(a\circ b),$$from which $u\in 
v\circ w$ for some $v\in a\circ x$, $w\in a\circ b$. We have $u\in 
v\circ w\subseteq \{v\}*(a\circ b)$ and $f{\Big(}\{v\}*(a\circ 
b){\Big)}\ge f(a)$, thus we have $f(u)\ge f(a)$, and $f$ is a fuzzy 
right ideal of $H$. We also have $f(u)\ge f(b)$. Indeed: Since $b\in 
H$ and $H$ is regular, there exists $y\in H$ such that 
$b\in\{b\}*\{y\}*\{b\}$. Then we have$$u\in a\circ b\subseteq 
\{a\}*\{b\}*\{y\}*\{b\}=(a\circ b)*(y\circ b).$$Then $u\in s\circ t$ 
for some $s\in a\circ b$, $t\in y\circ b$. Then we have$$u\in s\circ 
t\subseteq (a\circ b)*\{t\}=\{a\}*(b\circ t).$$ Since 
$f{\Big(}\{a\}*(b\circ t){\Big)}\ge f(b)$, we obtain $f(u)\ge f(b)$, 
and $f$ is a fuzzy left ideal of $H$. $\hfill\Box$\\From the 
Propositions 2.4 and 2.8 we have the following\medskip

\noindent{\bf Theorem 2.9.} {\it In regular hypersemigroups the 
concepts of fuzzy ideals and of fuzzy interior ideals 
coincide.}\medskip

\noindent{\bf Definition 2.10.} An hypersemigroup $H$ is called {\it 
intra-regular} if for every $a\in H$ there exist $x,y\in H$ such that 
$a\in (x\circ a)*(a\circ y)$.\medskip

\noindent{\bf Proposition 2.11.} {\it Let H be an hypersemigroup. The 
following are equivalent:

$(1)$ H is intra-regular.

$(2)$ $a\in H*\{a\}*\{a\}*H$ for every $a\in H$.

$(3)$ $A\subseteq H*A*A*H$ for every nonempty subset of 
$H$.}\medskip

\noindent{\bf Proposition 2.12.} {\it Let H be an intra-regular 
hypersemigroup and A an interior ideal of H. Then A is a subsemigroup 
of H.}\medskip

\noindent{\bf Proof.} Since $A$ is an interior ideal of $H$, we have 
$H*A*H\subseteq A$. Since $H$ is intra-regular, we have $A\subseteq 
H*A*A*H$. Then we have\begin{eqnarray*}A*A&\subseteq&(H*A*A*H)*A
=(H*A)*A*(H*A)\\&\subseteq&(H*H)*A*(H*H)\\&\subseteq&H*A*H\subseteq 
A.\end{eqnarray*}$\hfill\Box$\\By Propositions 2.7 and 2.12, we have 
the following\medskip

\noindent{\bf Corollary 2.13.} {\it In regular and in intra-regular 
hypersemigroups the interior ideals and the subidempotent interior 
ideals coincide.} \medskip

\noindent{\bf Proposition 2.14.} {\it Let H be an intra-regular 
hypersemigroup and f is a fuzzy interior ideal of H. Then f is a 
fuzzy ideal of H}.\medskip

\noindent{\bf Proof.} Let $a,b\in H$ and $u\in a\circ b$. Since $a\in 
H$ and $H$ is intra-regular, there exist $x,y\in H$ such that 
$a\in\{x\}*\{a\}*\{a\}*\{y\}$. Then$$a\circ b\subseteq 
\{x\}*\{a\}*\{a\}*\{y\}*\{b\}=(x\circ a)*{\Big(}(a\circ 
y)*\{b\}{\Big)}.$$Then $u\in v\circ w$ for some $v\in x\circ a$, 
$w\in (a\circ y)*\{b\}$. We have$$u\in v\circ w\subseteq (x\circ 
a)*\{w\}.$$Since $f$ is a fuzzy interior ideal of $H$, 
$f{\Big(}(x\circ a)*\{w\}{\Big)}\ge f(a)$. Thus we get $f(u)\ge 
f(a)$, and $f$ is a fuzzy right ideal of $H$. Since $b\in H$ and $H$ 
is intra-regular, there exist $z,t\in H$ such that $b\in 
\{z\}*\{b\}*\{b\}*\{t\}$, then we have$$a\circ b\subseteq 
\{a\}*\{z\}*\{b\}*\{b\}*\{t\}={\Big(}(a\circ z)*\{b\}{\Big)}*(b\circ 
t).$$Then $u\in c\circ d$ for some $c\in (a\circ z)*\{b\}$, $d\in 
b\circ t$. Then $u\in c\circ d\subseteq \{c\}*(b\circ t)$ and 
$f{\Big(}\{c\}*(b\circ t){\Big)}\ge f(b)$, so $f(u)\ge f(b)$, and $f$ 
is a fuzzy right ideal of $H$. $\hfill\Box$\\By Propositions 2.4 and 
2.14, we have the following theorem\medskip

\noindent{\bf Theorem 2.15.} {\it In intra-regular hypersemigroups 
the concepts of fuzzy ideals and fuzzy interior ideals 
coincide}.\medskip

\noindent{\bf Definition 2.16.} An ideal $A$ of an hypergroupoid $H$ 
is called {\it proper} if $A\not=H$. An hypergroupoid $H$ is called 
{\it simple} if does not contain proper ideals, that is, for every 
ideal $A$ of $H$, we have $A=H$.\medskip

\noindent{\bf Lemma 2.17.} {\it An hypersemigroup H is simple if and 
only if, for every $a\in H$, we have$$H=H*\{a\}*H.$$}{\bf Proof.} 
$\Longrightarrow$. Let $a\in H$. The set $H*\{a\}*H$ is an ideal of 
$H$. Indeed, it is a nonempty subset of $H$, and we have

$H*(H*\{a\}*H)=(H*H)*\{a\}*H\subseteq H*\{a\}*H$ and

$(H*\{a\}*H)*H=H*\{a\}*(H*H)\subseteq H*\{a\}*H$.\\
Since $H$ is simple, we have $H*\{a\}*H=H$.\\$\Longleftarrow$. Let 
$A$ be an ideal of $H$. Then $A=H$. Indeed: Take an element $x\in A$ 
$(A\not=\emptyset)$. By hypothesis, we have $H=H*\{x\}*H\subseteq 
(H*A)*H\subseteq A*H\subseteq A,$ thus we have $H=A$. 
$\hfill\Box$\medskip

\noindent{\bf Definition 2.18.} (cf. also [2]) An hypergroupoid $H$ 
is called {\it fuzzy simple} if every fuzzy ideal of $H$ is a 
constant function, that is, for every fuzzy ideal $f$ of $H$ and 
every $a,b\in H$, we have $f(a)=f(b)$.\medskip

\noindent{\bf Notation 2.19.} Let $H$ be an hypergroupoid and $a\in 
H$. We denote by $I_a$ the subset of $H$ defined as 
follows:$$I_a=\{b\in H \mid f(b)\ge f(a)\}.$${\bf Proposition 2.20.} 
{\it Let H be an hypergroupoid and f a fuzzy right ideal of H. Then 
the set $I_a$ is a right ideal of H for every $a\in H$.}\medskip

\noindent{\bf Proof.} Let $a\in H$. The set $I_a$ is a right ideal of 
$H$. Indeed: Since $a\in I_a$, the set $I_a$ is a nonempty subset of 
$H$. Moreover, $I_a*H\subseteq I_a$. Indeed: Let $x\in I_a*H$. Then 
$x\in H$ and $x\in u\circ v$ for some $u\in I_a$, $v\in H$. Since 
$x\in u\circ v$ and $f$ is a fuzzy right ideal of $H$, we have 
$f(x)\ge f(u)$. Since $u\in I_a$, we have $f(u)\ge f(a)$. Then we 
have $f(x)\ge f(a)$. Since $x\in H$ and $f(x)\ge f(a)$, we have $x\in 
I_a$. Thus $I_a$ is a right ideal of $H$. $\hfill\Box$\\
In a similar way we prove the following\medskip

\noindent{\bf Proposition 2.21.} {\it Let H be an hypergroupoid and f 
a fuzzy left ideal of H. Then the set $I_a$ is a left ideal of H for 
every $a\in H$.}\\
By Propositions 2.20 and 2.21 we have the following 
proposition\medskip

\noindent{\bf Proposition 2.22.} {\it Let H be an hypergroupoid and f 
a fuzzy ideal of H. Then the set $I_a$ is an ideal of H for every 
$a\in H$.}\medskip

\noindent{\bf Theorem 2.23.} {\it An hypergroupoid H is simple if and 
only if it is fuzzy simple.}\medskip

\noindent{\bf Proof.} $\Longrightarrow$. Let $f$ be a fuzzy ideal of 
$H$ and $a,b\in H$. Then $f(a)=f(b)$. Indeed: Since $f$ is a fuzzy 
ideal of $H$ and $a\in H$, by Proposition 2.22, the set $I_a=\{b\in H 
\mid f(b)\ge f(a)\}$ is an ideal of $H$. Since $H$ is simple, we have 
$I_a=H$. Then $b\in I_a$, so $f(b)\ge f(a)$. By symmetry, we get 
$f(a)\ge f(b)$, thus $f(a)=f(b)$.\\
$\Longleftarrow$. Let $H$ be fuzzy simple and $I$ an ideal of $H$. 
Then $I=H$. Indeed: Let $x\in H$. Since $I$ is an ideal of $H$, the 
characteristic function $f_I$ is a fuzzy ideal of $H$. Since $H$ is
fuzzy simple, $f_I$ is a constant function, that is, $f_I(y)=f_I(z)$ 
for every $y,z\in H$. Take an element $a\in I$ $(I\not=\emptyset)$. 
Then we have $f_I(x)=f_I(a)=1$, so $x\in I$. Thus $H$ is simple. 
$\hfill\Box$\medskip

\noindent{\bf Theorem 2.24.} {\it Let H be an hypersemigroup. Then H 
is simple if and only if every fuzzy interior ideal of H is a 
constant function.}\medskip

\noindent{\bf Proof.} $\Longrightarrow$. Let $f$ be a fuzzy interior 
ideal of $H$ and $a,b\in H$. Then $f(a)=f(b)$. Indeed: Since $a\in H$ 
and $H$ is simple, by Lemma 2.17, we have $H=H*\{a\}*H$. Since $b\in 
H$, we have $b\in (x\circ a)*\{y\}$ for some $x,y\in H$. Since $f$ is 
a fuzzy interior ideal of $H$, we have $f(b)\ge f(a)$. By symmetry, 
we get $f(a)\ge f(b)$, so $f(a)=f(b)$.\\
$\Longleftarrow$. Let $f$ is a fuzzy ideal of $H$. By Proposition 
2.4, $f$ is a fuzzy interior ideal of $H$. By hypothesis, $f$ is a 
constant function. Thus $H$ is fuzzy simple. Then, by Theorem 2.23, 
$H$ is simple. $\hfill\Box$\\
As a consequence we have the following:\medskip

\noindent{\bf Theorem 2.25.} {\it If H is an hypersemigroup H, then 
the following are equivalent:

$(1)$ H is simple.

$(2)$ $H=H*\{a\}*H$.

$(3)$ H is fuzzy simple.

$(4)$ Every fuzzy interior ideal of $H$ is a constant 
function.}{\small\bigskip

\bigskip

\noindent University of Athens\\
Department of Mathematics\\
15784 Panepistimiopolis, Greece\\
email: nkehayop@math.uoa.gr

\end{document}